\documentclass[a4paper,12pt,number]{elsart}

\usepackage[utf8]{inputenc}
\usepackage{amsmath}
\usepackage{amsfonts}
\usepackage{algorithm}
\usepackage{algpseudocode} 
\usepackage{graphicx}
\usepackage{url}

\begin{document}
\begin{frontmatter}

\title{Accurate Bidiagonal Decomposition and Computations with Generalized Pascal Matrices}

\author[Zaragoza] {J. Delgado \thanksref{BFM}}
\author[Zaragoza] {H. Orera \thanksref{BFM}}
\author[Zaragoza] {J. M. Pe\~na \thanksref{BFM}}

\address[Zaragoza] {Departamento de Matem\'{a}tica Aplicada, Universidad de Zaragoza, Spain}

\thanks[BFM] {This work was partially supported through the Spanish research grant PGC2018-096321-B-I00
(MCIU/AEI), by Gobierno de Arag\'on (E41-17R) and Feder 2014-2020 ``Construyendo Europa desde
Arag\'on''.}
\begin{abstract}
This paper provides an accurate method to obtain the bidiagonal factorization of many generalized Pascal matrices, which in turn can be used to compute with high relative accuracy the eigenvalues, singular values and inverses of these matrices. Numerical examples are included.


\end{abstract}
\end{frontmatter}

\section{Introduction}

Finding classes of matrices relevant in applications for which algebraic computations can be performed with high relative accuracy (HRA) is an active research topic of great interest in recent years. This goal has been achieved for some subclasses of totally positive matrices (see, for instance,  \cite{DK,  SBV, BV2, JS, qBer, dop1, dop2}). Let us recall that a matrix is called {\it totally positive} (TP) if all their minors are nonnegative and, if they are all positive, then the matrix is called {\it strictly totally positive} (STP). These classes of matrices play an important role in many fields such as approximation theory, statistics, mechanics, computer-aided geometric design, economics, combinatorics or biology (see \cite{A, GM, p}). For the subclasses of TP matrices mentioned above, their bidiagonal decomposition (see Section 2) was obtained with HRA, and then the algorithms given in \cite{K1, K2, KoevSoft} permit to compute with HRA many algebraic calculations: all their eigenvalues and singular values, their inverses, or the solution of some linear systems. Recall that a  real value $x$ is obtained with HRA if the relative error
of the computed value  $\tilde x$ satisfies
$ 
 \| x-\tilde x \|/  \| x \| < Ku,
$
where $K$ is a positive constant independent of the arithmetic precision and $u$ is the unit round-off. It is well known that an algorithm can be performed with HRA if all the included subtractions are of initial data, that is, if it only includes products, divisions, sums of numbers of the same sign and subtractions of the initial data (cf. \cite{DK, AN, K2}).

The lower triangular Pascal matrix $P=(p_{ij})_{1\leq i,j \leq n+1}$ (with $p_{ij}=\binom{i-1}{j-1}$ for $1\leq j \leq i \leq n+1$ and $p_{ij}:=0$ when $j>i$) and the symmetric Pascal matrix $R=(r_{ij})_{1\leq i,j \leq n+1}$ (with $r_{ij}=\binom{i+j-2}{j-1}$) are naturally derived from the Pascal
triangle. The matrix $R=PP^T$ is also called Pascal matrix. This paper deals with some classes of matrices (see \cite{Ref1, Ref2, Ref3, Ref5}) generalizing  the lower triangular Pascal matrix and the symmetric Pascal matrix. These generalized classes of Pascal matrices arise in applications such as filter design, probability, combinatorics, signal processing or electrical engineering (see \cite{X} and its references). In \cite{X}, one can also see some concrete applications of solving linear systems with these matrices.
Let us recall that the bidiagonal decomposition of a Pascal matrix is well known and has the remarkable property that it is formed by 1's (see \cite{K2, a}). In this paper, we show that the bidiagonal decompositions of these generalized Pascal matrices can be obtained with HRA, and so the remaining algebraic calculations mentioned above can be also computed with HRA. Although Pascal matrices are ill-conditioned (see \cite{a}) and the bidiagonal decompositions of their generalizations are not as simple as those with the Pascal matrix, we can still guarantee the mentioned algebraic calculations with HRA. 

In Section 2 we present auxiliary results concerning the bidiagonal decomposition of nonsingular TP matrices and some basic definitions of generalized Pascal matrices. In Section 3 we obtain the bidiagonal decomposition of generalized triangular Pascal matrices and of lattice path matrices, which in turn contain many classical generalized Pascal matrices. In many cases,  we prove that they are TP or STP and show that the algebraic calculations mentioned above can be computed with HRA. Section 4 includes numerical experiments showing the great accuracy of the proposed method. Finally, Section 5 summarized the main conclusions of the paper.

\section{Auxiliary results and basic definitions}
\emph{Neville elimination} (NE) is an alternative procedure to Gaussian elimination that produces zeros in a column of a matrix by adding to each row an appropriate multiple of the previous one. This elimination procedure is very useful when dealing with some classes of matrices such as TP matrices. For more details on NE see \cite{TPandNE,fact}. Given a nonsingular matrix $A=(a_{ij})_{1\leq i,j \leq n}$, the Neville elimination procedure consists of $n-1$ steps and leads to the following sequence of matrices:
\begin{equation} \label{seq}
A=:A^{(1)}\rightarrow \widetilde{A}^{(1)}\rightarrow A^{(2)}\rightarrow \widetilde{A}^{(2)}
        \rightarrow \cdots \rightarrow A^{(n)}= \widetilde{A}^{(n)}=U,
\end{equation}
where $U$ is an upper triangular matrix.

The matrix $\widetilde{A}^{(k)}=(\widetilde{a}_{ij}^{(k)})_{1\leq i,j \leq n}$ is obtained from the matrix $A^{(k)}=(a_{ij}^{(k)})_{1\leq i,j \leq n}$  by a row permutation that moves to the bottom the rows with a zero entry  in column $k$ below the main diagonal. For nonsingular TP matrices, it is always possible to perform NE without row exchanges (see \cite{TPandNE}). If a row permutation is not necessary at the $k$th step, we have that $\widetilde{A}^{(k)}=A^{(k)}$.
 The entries of  $A^{(k+1)}=(a_{ij}^{(k+1)})_{1\leq i,j \leq n}$ can be obtained 
from $\widetilde{A}^{(k)}=(\widetilde{a}_{ij}^{(k)})_{1\leq i,j \leq n}$ using the formula:

\begin{equation}
a_{ij}^{(k+1)} = \left\{ \begin{array}{ll}
    \widetilde{a}_{ij}^{(k)}   -  \dfrac{ \widetilde{a}_{ik}^{(k)} }{\widetilde{a}_{i-1,k}^{(k)}}\widetilde{a}_{i-1,j}^{(k)},  & \text{if} \ k  \leq j<i \leq n  \text{ and }  \widetilde{a}_{i-1,k}^{(k)} \not= 0,\\
    \widetilde{a}_{ij}^{(k)},   & \text{otherwise,}
   \end{array}\right.
\end{equation}
for $k=1,\ldots,n-1$. The $(i,j)$ \textit{pivot} of the NE of $A$ is given by
 \[ p_{ij} = \widetilde{a}_{ij}^{(j)}, \ \  1 \leq j \leq i \leq n. \]
If $i=j$  we say that $p_{ii}$ is a \textit{diagonal pivot}. The $(i,j)$ \emph{multiplier}  of the NE of $A$, with $1\leq j\leq i \leq n$,  is defined as

\begin{equation*}
m_{ij}=\left\{ \begin{array}{ll}
   \dfrac{\widetilde{a}_{ij}^{(j)}}{\widetilde{a}_{i-1,j}^{(j)}}=\dfrac{p_{ij}}{p_{i-1,j}},  & \text{if}  \ \widetilde{a}_{i-1,j}^{(j)}\not= 0, \\
  0,  & \text{if}  \  \widetilde{a}_{i-1,j}^{(j)} = 0.
   \end{array}\right.
\end{equation*}
The multipliers  satisfy  that
\[m_{ij} = 0\Rightarrow m_{hj} = 0 \ \ \ \forall h > i.\]
Nonsingular TP matrices can be expressed as a product of nonnegative bidiagonal matrices.  The following theorem (see Theorem 4.2 and p. 120 of \cite{fact})  introduces this representation, which is called  the \textit{bidiagonal decomposition}.
\begin{thm}(cf. Theorem 4.2  of \cite{fact})\label{BD} Let $A=(a_{ij})_{1\leq i,j\leq n}$  be a nonsingular TP matrix. Then $A$ admits the following representation:
\begin{equation}\label{BD.definicion}
A=F_{n-1}F_{n-2}\cdot \cdot \cdot F_{1}DG_{1}\cdot \cdot \cdot G_{n-2}G_{n-1}, 
\end{equation}
where $D$ is the diagonal matrix \emph{diag}$(p_{11},\ldots,p_{nn})$ with positive diagonal entries  and $F_i$, $G_i$ are the nonnegative bidiagonal matrices given by

\begin{equation}\label{BD.F}
F_{i}=\left(\begin{matrix}
1 &        &   &  &  &   &  \\
0 & 1      &   &  &  &   & \\
  & \ddots & \ddots  &  &   &  & \\ 
 &  & 0 & 1 &   &  & \\
 & &  & m_{i+1,1} &  1 &  & \\
 & &  &  & \ddots  & \ddots & \\
 & &  &  &   & m_{n,n-i} & 1
\end{matrix}\right),
\end{equation}
\\
\begin{equation}\label{BD.G}
G_{i}=\left(\begin{matrix}
1 &   0     &   &  &  &   &  \\
 & 1      &  \ddots &  &  &   & \\
  &  & \ddots  & 0 &   &  & \\ 
 &  &  & 1 &  \widetilde{m}_{i+1,1} &  & \\
 & &  & &  1 & \ddots  & \\
 & &  &  &   & \ddots & \widetilde{m}_{n,n-i}\\
 & &  &  &   &  & 1
\end{matrix}\right),
\end{equation}
\\
for all $i \in \{ 1,\ldots,n - 1 \}$. If, in addition,  the entries $m_{ij}$ and $\widetilde{m}_{ij}$ satisfy 
\begin{equation}\label{BD.unique}
\begin{array}{c}
m_{ij} = 0\Rightarrow m_{hj} = 0 \quad \forall h > i,\\
\widetilde{m}_{ij} = 0\Rightarrow \widetilde{m}_{hj} = 0 \quad\forall h > i,
\end{array}
\end{equation}
then the decomposition is unique.
\end{thm}
In the bidiagonal decomposition given by \eqref{BD.definicion}, \eqref{BD.F} and \eqref{BD.G}, the entries $m_{ij}$ and $p_{ii}$ are  the multipliers and diagonal pivots, respectively,  corresponding to the NE of $A$ (see Theorem 4.2  of \cite{fact} and the comment below it)  and the entries $\widetilde{m}_{ij}$ are the multipliers of the NE of $A^T$ (see p. 116 of \cite{fact}). The following result shows that the bidiagonal decomposition also characterizes STP matrices. 
\begin{thm}\label{char.STP} (cf. Theorem 4.3  of \cite{fact})
A nonsingular $n\times n$ matrix $A$ is STP if and only if it can be factorized in the form \eqref{BD.definicion} with $D$ a diagonal matrix with positive diagonal entries, $F_i$, $G_i$ given by \eqref{BD.F} and \eqref{BD.G}, and the entries $m_{ij}$ and $\widetilde{m}_{ij}$ are positive numbers. This factorization is unique.
\end{thm}

 Let us recall that an algorithm can be performed with high relative accuracy if it only includes products, divisions, sums of numbers of the same sign and subtractions of initial data (cf. \cite{AN,DK}). In \cite{K1,K2}, assuming that the bidiagonal decomposition of a nonsingular TP matrix $A$ is known to HRA, Plamen Koev designed efficient algorithms for computing to HRA the eigenvalues, singular values and the inverse of $A$ as well as the solution to linear systems of equations $Ax=b$ whenever $b$ has a pattern of alternating signs.

In \cite{K1} the matrix notation $\mathcal{BD}(A)$ was introduced to represent the bidiagonal decomposition of a nonsingular TP matrix,
 \begin{equation}\label{eq:BD.A}
 	(\mathcal{BD}(A))_{ij}=\left\{
 		\begin{array}{ll}
 			m_{ij}, & \text{if } i>j, \\
 			\widetilde{m}_{ji}, & \text{if } i<j, \\
 			p_{ii}, & \text{if } i=j.
 		\end{array}
	\right.
 \end{equation}

In general, more matrices can be written as a product of bidiagonal matrices following \eqref{BD.definicion}. Throughout this paper, we will use the notation $\mathcal{BD}(A)$ given by \eqref{eq:BD.A} to denote the bidiagonal decomposition \eqref{BD.definicion}-\eqref{BD.G} of a general matrix $A$.

  Finally, let us introduce the following classical generalizations of Pascal matrices. 
  \begin{defn}\label{def1}(see \cite{Ref1,Ref5})
  For a real number $x$, the generalized Pascal matrix of the first kind, $P_n[x]$, is defined as the $(n+1)\times (n+1)$ lower triangular matrix with $1'$s on the main diagonal and
\[
(P_n[x])_{ij}:=x^{i-j}\binom{i-1}{j-1}, \quad 1\leq j \leq i \leq n+1 
\]
and the symmetric generalized Pascal $(n+1)\times (n+1)$ matrix  $R_n[x]$ is given by
\[(R_n[x])_{ij}:=x^{i+j-2}\binom{i+j-2}{j-1}, \quad 1\leq i,j  \leq n+1. \]
For $x,y\in\mathbb{R}$ we define the $(n+1)\times (n+1)$ matrix $R_n[x,y]$
\[(R_n[x,y])_{ij}:=x^{j-1}y^{i-1}\binom{i+j-2}{j-1}, \quad 1\leq i,j  \leq n+1. \]
  \end{defn}
 Observe that $R_n[x]=R_n[x,x]$, that $P_n[1]$ is the lower triangular Pascal matrix and that $R_n[1]$ is the symmetric Pascal matrix.
 
  \begin{defn}\label{def2}(see \cite{Ref2}) For $x,y\in\mathbb{R}$, the extended generalized Pascal matrix $\Phi_n[x, y]$ is defined as
\[  
(\Phi_n[x,y])_{ij}=x^{i-j}y^{i+j-2}\binom{i-1}{j-1}, \quad 1\leq j \leq i \leq n+1  \] 
and the extended generalized symmetric Pascal matrix $\Psi_n[x,y]$ is given by
\[  (\Psi_n[x, y])_{ij}=x^{i-j}y^{i+j-2}\binom{i+j-2}{j-1}, \quad 1\leq i,j  \leq n+1. \]
 \end{defn}
In the next section we are going to deduce the bidiagonal decomposition of more general classes of matrices. As a consequence, we can also obtain the bidiagonal decomposition of the matrices $P_n[x]$, $R_n[x,y]$, $\Phi_n[x,y]$ and $\Psi_n[x,y]$.

\section{Bidiagonal decomposition of generalized Pascal matrices}

\subsection{Generalized triangular Pascal matrices}

 Let $x$ and $\lambda$ be two real numbers and let $n$ be a nonnegative integer. We define the notation $x^{n|\lambda}$ as follows:
\begin{equation}\label{notacion}
x^{n|\lambda}=\left\{\begin{matrix}
x(x + \lambda)\cdots (x + (n-1)\lambda), & \text{ if } n>0,\\ 
1, & \text{ if } n=0.
\end{matrix}\right. 
\end{equation}
 In \cite{Ref3}, the generalized lower triangular Pascal matrix $P_{n,\lambda}[x]$ is defined by
 \begin{equation}\label{gen3}
 (P_{n,\lambda}[x])_{i,j}:=x^{(i-j)|\lambda}\binom{i-1}{j-1}, \quad 1\leq j\leq i \leq n+1,
 \end{equation}
where $n$ is a natural number and $\lambda$ and $x$ are  real numbers. Observe that the particular case $\lambda=0$ leads to the generalized Pascal matrix of the first kind $P_{n,0}[x]=P_n[x]$. The following result provides the bidiagonal decomposition of the generalized Pascal matrix $P_{n,\lambda}[x]$.

 \begin{thm}\label{BD3} Given $x,\lambda\in\mathbb{R}$ and $n\in \mathbb{N}$,
 let $P_{n,\lambda}[x]$ be the $(n+1) \times (n+1)$ lower triangular matrix given by \eqref{gen3}. 
 
  \begin{itemize}
  \item[i)] If $x\neq k \lambda$ for $k=-n+1,\ldots,0,\ldots,n-1$, then we have that
 \begin{equation}\label{BDgen3}
 \left(\mathcal{BD}(P_{n,\lambda}[x] )\right)_{ij}=\left\{\begin{matrix} 
 1, &  i= j ,\\ 
 x+(i-2j)\lambda, &  i> j ,\\ 
0, &  i<j.
\end{matrix}\right.
 \end{equation}
 \item[ii)] If $x= k \lambda$ for some $k\in\{0,\ldots,n-1\}$, then we have that
 \begin{equation}\label{BDgen3.1}
\left(\mathcal{BD}(P_{n,\lambda}[x] )\right)_{ij}=\left\{\begin{matrix} 
 1, &  i= j ,\\ 
 x+(i-2j)\lambda, &  i> j, j\leq k,\\ 
0, &  \emph{otherwise}.
\end{matrix}\right.
\end{equation}
\item[iii)] If $x=- k \lambda$ for some $k\in\{0,\ldots,n-1\}$, then we have that
 \begin{equation}\label{BDgen3.2}
\left(\mathcal{BD}(P_{n,\lambda}[x] )\right)_{ij}=\left\{\begin{matrix} 
 1, &  i= j ,\\ 
 x+(i-2j)\lambda, &  0 <i-j\leq k,\\ 
0, &  \emph{otherwise}.
\end{matrix}\right.
 \end{equation}
 \end{itemize}
  
 \end{thm}

 \begin{pf}
Let us first assume that  $x\neq k \lambda$ for $k=-n+1,\ldots,0,\ldots,n-1$. We are going to perform the first step of the Neville elimination of $A=(a_{ij})_{1\leq i,j \leq n+1}$, where $a_{ij}:=(P_{n,\lambda}[x])_{i,j}$ for $i,j=1,\ldots,n+1$:
 \[a_{ij}^{(2)}=a_{ij}-\frac{a_{i1}}{a_{i-1,1}}a_{i-1,j}=a_{ij}-(x+(i-2)\lambda)a_{i-1,j}, \quad i>j\geq 1.\]
Applying \eqref{gen3} to the previous formula, $a_{ij}^{(2)}$ can be written as
\begin{equation*}
a_{ij}^{(2)}=x^{(i-j)|\lambda}\binom{i-1}{j-1}- (x+(i-2)\lambda)x^{(i-j-1)|\lambda}\binom{i-2}{j-1}.
\end{equation*}
By formula $\eqref{notacion}$, we have that
 \begin{align*}
a_{ij}^{(2)}&=\left((x+(i-j-1)\lambda)\binom{i-1}{j-1}- (x+(i-2)\lambda)\binom{i-2}{j-1} \right)x^{(i-j-1)|\lambda}\\
&=\left(x\binom{i-2}{j-2}+\frac{(i-j-1)(i-1)!}{(j-1)!(i-j)!}\lambda - \frac{(i-2)(i-2)!}{(j-1)!(i-j-1)!}\lambda \right)x^{(i-j-1)|\lambda}.\\
\end{align*}
After some computations we deduce that
\begin{equation*}
a_{ij}^{(2)}=\left(x\binom{i-2}{j-2}-\lambda\binom{i-2}{j-2}\right)x^{(i-j-1)|\lambda}=\binom{i-2}{j-2}(x-\lambda)^{(i-j)|\lambda}.
\end{equation*}
We can observe that $a_{ij}^{(2)}=(P_{n,\lambda}[x])_{ij}^{(2)}=(P_{n,\lambda}[x-\lambda])_{i-1,j-1}$ for $i>j\geq 2$ and, hence, $(P_{n,\lambda}[x])^{(2)}[2,\ldots,n+1]=(P_{n,\lambda}[x-\lambda])[1,\ldots,n]$. Then we can deduce that  $(P_{n,\lambda}[x])_{ij}^{(k+1)}=(P_{n,\lambda}[x-k\lambda])_{i-k,j-k}$ for $i>j\geq k+1$ and that the multipliers for the $k$th step of the NE of $P_{n,\lambda}[x]$ will be given by $x-(k-1)\lambda+(i-k-1)\lambda$ for $i=k+1,\ldots,n+1$, and so, we conclude that \eqref{BDgen3} holds. 

Let us now assume that $x=k\lambda$ for any $k\in\{0,\ldots,n-1\}$. Following the above proof we can see that $(P_{n,\lambda}[x])_{ij}^{(k+1)}=(P_{n,\lambda}[0])_{i-k,j-k}$ and the NE finishes at the $k+1$ step. Hence, $ii)$ holds.

Finally, if $x=-k\lambda$ for any $k\in\{0,\ldots,n-1\}$, then $x^{(i-j)|\lambda}=0$ for $i-j>k$. Then the $n-k$ lower subdiagonals are already zero and the associated multipliers are also zero since the elimination procedure is not carried out on those entries. So, we conclude that $iii)$ holds.
 \end{pf}
 
 \begin{rem}\label{rem1}
It can be checked that the computational cost for the bidiagonal decomposition in \eqref{BDgen3} is of $\mathcal{O}\left(n^2\right)$ elementary operations. For the bidiagonal decompositions in \eqref{BDgen3.1} and \eqref{BDgen3.2} the computational costs are of $\mathcal{O}\left(k^2\right)$ and of  $\mathcal{O}\left(k\cdot n\right)$ elementary operations, respectively.
 \end{rem}

 The following corollary characterizes the matrices $P_{n,\lambda}[x]$ that are TP.
 
 \begin{cor}\label{GPT.TP}
Let $P_{n,\lambda}[x]$ be the lower triangular matrix given by \eqref{gen3}  with $x,\lambda\in\mathbb{R}$ and with $n\in \mathbb{N}$.  Then $P_{n,\lambda}[x]$ is a TP matrix if and only if one of the following conditions holds:
\begin{itemize}
\item[i)] $x\geq (n-1)|\lambda|$.
\item[ii)] $x=k|\lambda|$ for $k=0,\ldots,n-1$.
\end{itemize}
 \end{cor}
 \begin{pf}
 By Theorem \ref{BD3} we know that $P_{n,\lambda}[x]$ admits a factorization as a product of bidiagonal matrices. If $i)$ or $ii)$ holds, then all the bidiagonal matrices are nonnegative and so TP, and hence, its product is also TP (see for example Theorem 3.1 of \cite{A}). Conversely, if $P_{n,\lambda}[x]$ is TP, 
since it is also nonsingular, it admits a unique bidiagonal decomposition by Theorem \ref{BD}. Moreover, this bidiagonal decomposition will be given by Theorem \ref{BD3} and the $m_{ij}$'s will be nonnegative. Hence, either $i)$ or $ii)$ holds.
 \end{pf}
  The previous definition of $P_{n,\lambda}[x]$ is generalized in \cite{Ref3} for two variables $x,y$ as follows:
 \begin{equation}\label{gen3b}
 (P_{n,\lambda}[x,y])_{i,j}:=x^{(i-j)|\lambda}y^{(j-1)|\lambda}\binom{i-1}{j-1}.
 \end{equation}
 Let us also define $P_n[x,y]:=P_{n,0}[x,y]$. It is straightforward to see that the matrix $P_{n,\lambda}[x,y]$ can be expressed as the product of $P_{n,\lambda}[x]$ and a diagonal matrix:
 \begin{equation}\label{diag1}
  P_{n,\lambda}[x,y] = P_{n,\lambda}[x] \ \text{diag}(1,y^{1|\lambda},\ldots,y^{n|\lambda}) .
 \end{equation}
 In \cite{Ref3} a further generalization of $P_{n,\lambda}[x,y]$ is given in terms of an arbitrary sequence $\textbf{a}=\{a_n\}_{n\geq 0}$
\begin{equation}\label{gen3c}
 (P_{n,\lambda}[x,y,\textbf{a}])_{i,j}:=a_{j-1}x^{(i-j)|\lambda}y^{(j-1)|\lambda}\binom{i-1}{j-1},
 \end{equation} 
 and so we also derive 
  \begin{equation}\label{diag2}
  P_{n,\lambda}[x,y,\textbf{a}] = P_{n,\lambda}[x] \ \text{diag}(a_0,a_1 y^{1|\lambda},\ldots,a_n y^{n|\lambda}) .
 \end{equation}
 Observe that the matrix $P_{n,\lambda}[x,y] =P_{n,\lambda}[x,y,\textbf{1}]$, where $\textbf{1}$ is the sequence formed by $1's$.
 By \eqref{diag2} and Theorem $\ref{BD3}$, we can deduce the bidiagonal decomposition of the matrix $\mathcal{BD}(P_{n,\lambda}[x,y,\textbf{a}])$. For example, if $x\neq k\lambda$ for $k=-n+1,\ldots,0,\ldots,n-1$, its bidiagonal decomposition is given by
  \begin{equation}\label{BDgen3a}
 \left(\mathcal{BD}(P_{n,\lambda}[x,y,\textbf{a}] )\right)_{ij}=\left\{\begin{matrix} 
 a_{j-1}y^{(j-1)|\lambda}, &  i= j ,\\ 
 x+(i-2j)\lambda, &  i> j ,\\ 
0, &  i<j.
\end{matrix}\right.
 \end{equation}
 
 \subsection{Lattice path matrices}
 
 Let  $Lp_n(\alpha,\beta,\gamma)=(k_{ij})_{1\leq i,j \leq n+1}$ be the $(n+1)\times (n+1)$ {\it lattice path matrix} 
 such that its entries are given by the recurrence relation
 \begin{equation}\label{rec}
  \alpha k_{i,j-1} + \beta k_{i-1,j} + \gamma k_{i-1,j-1}=k_{ij},\quad 2\leq i,j\leq n+1,
 \end{equation}
 with $k_{1j}=\alpha^{j-1}$ for $j\in\{1,\ldots,n+1\}$ and $k_{i1}=\beta^{i-1}$ for
 $i\in\{1,\ldots,n+1\}$. These matrices were considered in \cite{Ref5}. Other related classes of matrices were considered in \cite{B}, where it was also shown that some of those matrices are TP. In Theorem 2.3 of \cite{Ref5} it is also shown that $Lp_n(\alpha,\beta,\gamma)$ admits the following factorization
 \begin{equation}\label{factLDU}
 Lp_n(\alpha,\beta,\gamma)= P_n [\alpha] D^n_{\alpha\beta+\gamma}(P_n [\beta])^T,
 \end{equation}
 where $D^n_{\alpha\beta+\gamma}=\text{diag}(1,\alpha\beta+\gamma,\ldots,(\alpha\beta+\gamma)^n)$ and $P_n[\delta]=P_{n,0}[\delta]$.
Observe that the matrix $Lp_n(\alpha,\beta,\gamma)$ is nonsingular if and only if
$\alpha\,\beta+\gamma\ne 0$ In the following result, we deduce the bidiagonal 
decomposition of $Lp_n(\alpha,\beta,\gamma)$.
 
 \begin{thm}\label{BD.lattice}
Let $Lp_n(\alpha,\beta,\gamma)=(k_{ij})_{1\leq i,j \leq n+1}$ be the matrix whose entries are defined by \eqref{rec} with $\alpha\,\beta+\gamma\ne 0$. Then its bidiagonal decomposition is given by
\begin{equation}\label{BDgeneral}
\left(\mathcal{BD}( Lp_n(\alpha,\beta,\gamma))\right)_{ij}=\begin{cases}
(\alpha\beta+\gamma)^{i-1}, & \text{ if } i=j, \\ 
\alpha, & \text{ if } i>j, \\
\beta, & \text{ if } i<j.
\end{cases}
\end{equation}
 \end{thm}
 \begin{pf}
 By \eqref{factLDU}, the matrix $Lp_n(\alpha,\beta,\gamma)$ can be decomposed as the product of a lower triangular matrix, a diagonal matrix and an upper triangular matrix. Hence, we can deduce its bidiagonal decomposition from the bidiagonal decomposition of these three factors. Since $P_n[\alpha]=P_{n,0}[\alpha]$, by Theorem \ref{BD3} we have that
 \begin{equation*}
\left(\mathcal{BD}( P_n[\alpha])\right)_{ij}=\begin{cases}
1, & \text{ if } i=j, \\ 
\alpha, & \text{ if } i>j,\\
0, & \text{ if } i<j.
\end{cases}
\end{equation*}
Analogously, the bidiagonal decomposition of $( P_n[\beta])^T$ is given by
 \begin{equation*}
\left(\mathcal{BD}( (P_n[\beta])^T)\right)_{ij}=\begin{cases}
1, & \text{ if } i=j, \\ 
0, & \text{ if } i>j,\\
\beta, & \text{ if } i<j.
\end{cases}
\end{equation*}
Therefore, we conclude that  \begin{equation}\label{factLDU2}
 Lp_n(\alpha,\beta,\gamma)= P_n [\alpha] D^n_{\alpha\beta+\gamma}(P_n [\beta])^T=\bar F_{n}\cdots \bar F_1  D^n_{\alpha\beta+\gamma}\bar G_{1}\cdots \bar G_{n},
 \end{equation}
 where $\bar F_{k}$ is the lower bidiagonal matrix given by \eqref{BD.F} with all multipliers equal to $\alpha$  and $\bar G_{k}$ is the upper bidiagonal matrix given by \eqref{BD.G} with all multipliers equal to $\beta$. So, \eqref{BDgeneral} holds.
  \end{pf}
  The following corollary considers a case where $Lp_n(\alpha,\beta,\gamma)$ is STP and shows that its bidiagonal decomposition can be computed to HRA.
  \begin{cor}\label{corSTP}
Let $Lp_n(\alpha,\beta,\gamma)=(k_{ij})_{1\leq i,j \leq n+1}$ be the matrix whose entries are defined by \eqref{rec}. If $\alpha,\beta>0$ and $\alpha\beta+\gamma>0$, then $Lp_n(\alpha,\beta,\gamma)=(k_{ij})_{1\leq i,j \leq n+1}$ is an STP matrix. Moreover, if $\gamma\geq 0$, then its bidiagonal decomposition \eqref{BDgeneral} can be computed to HRA and it can be used to obtain the eigenvalues, singular values and the inverse of $Lp_n(\alpha,\beta,\gamma)$ with HRA as well as the solution of the linear systems $Lp_n(\alpha,\beta,\gamma)x=b$, where $b=(b_1,\ldots,b_{n+1})$ has alternating signs.
  \end{cor}
\begin{pf}
By Theorem \ref{char.STP}, $Lp_n(\alpha,\beta,\gamma)=(k_{ij})_{1\leq i,j \leq n+1}$ is an STP matrix. With the additional condition $\gamma\geq 0$, $\mathcal{BD}( Lp_n(\alpha,\beta,\gamma))$ can be computed with a subtraction-free algorithm, and hence, with HRA, which in turn guarantees that the algebraic computations stated in the statement of this corollary can be performed with HRA (see Section \ref{sec.num} or Section 3 of \cite{K2}). 
\end{pf}  
\begin{rem}\label{rem2}
In order to compute the bidiagonal decomposition \eqref{BDgeneral}, $n+1$ elementary operations are necessary, that is, a computational cost of $\mathcal{O}\left(n\right)$.
\end{rem}
 The class of lattice path matrices, $Lp_n(\alpha,\beta,\gamma)$, contains the generalizations of Pascal matrices given by definitions \ref{def1} and \ref{def2}, and so, from their bidiagonal decomposition  we can deduce the bidiagonal decomposition of these  matrices. In particular, in Theorem 3.1 of \cite{Ref5} the following relationship  was  proved:
\begin{equation}\label{casos}
Lp_n(\alpha,\beta,\gamma)=\left\{\begin{array}{ll}
P_n[x,y], & \text{ if } \alpha=0, \beta=y, \gamma=x,\\ 
R_n[x,y],&  \text{ if } \alpha=x, \beta=y, \gamma=0, \\ 
\Phi_n[x,y], & \text{ if } \alpha=0, \beta=xy, \gamma=y^2, \\  
\Psi_n[x,y], & \text{ if } \alpha=y/x, \beta=xy, \gamma=0. \\  
\end{array}\right.
\end{equation} 
Taking into account   \eqref{casos}, we can use  Theorem \ref{BD.lattice} to obtain their bidiagonal decomposition. We can also apply Corollary  \ref{corSTP} to study the cases when they are STP and when their  bidiagonal decomposition can be obtained with HRA.
\begin{cor}\hfill
\begin{itemize}
	\item[i)] If $x,y>0$, then $R_n[x,y]$ is STP and its bidiagonal decomposition can be computed to HRA.
	\item[ii)] If $xy>0$, then $\Psi_n[x,y]$ is STP and its bidiagonal decomposition can be computed to HRA.
\end{itemize}
\end{cor}
 
\section{Numerical experiments}\label{sec.num}

In  \cite{K1,K2}, assuming that the parameterization $\mathcal{BD}(A)$ of 
a nonsingular TP matrix $A$ is known, Plamen Koev 
presented algorithms to solve the following algebraic problems for $A$:
computation of the eigenvalues and the singular values of $A$, computation of 
$A^{-1}$ and solution of the systems of linear equations $Ax=b$. 
In \cite{MM.Inv} Marco and Mart\'inez presented another algorithm for the
computation of $A^{-1}$ from $\mathcal{BD}(A)$.
If, in addition, $\mathcal{BD}(A)$ is known to HRA, then the algorithms 
solve these algebraic problems to HRA (in the case of linear systems only
when $b$ has a pattern of alternating signs).
Koev implemented the corresponding algorithms for Matlab and Octave, which are available 
in the software library {\it TNTool} in \cite{KoevSoft}. 
The functions are \verb"TNEigenValues" for the eigenvalues, \verb"TNSingularValues"
for the singular values, \verb"TNInverseExpand" for the inverse and \verb"TNSolve"
for the solution of linear system of equations. 
The functions require as input argument the data determining 
the bidiagonal decomposition \eqref{BD.definicion}-\eqref{BD.G}  of $A$, or equivalently, 
$\mathcal{BD}(A)$ given by \eqref{eq:BD.A}, and, in the case of \verb"TNSolve",
in addition, the vector $b$.

\begin{rem}The computational cost for both \verb"TNSolve" and 
\verb"TNInverseExpand" is $\mathcal{O}\left(n^2\right)$ elementary operations (see \cite{K2} and Section 4 of \cite{MM.Inv}) and for
the other two functions, \verb"TNEigenValues" and \verb"TNSingularValues", is $\mathcal{O}\left(n^3\right)$
 elementary operations. Hence, taking into account
remarks \ref{rem1} and \ref{rem2}, the total computational cost for solving a linear system or computing the inverse with the
matrices corresponding to these bidiagonal computations is $\mathcal{O}\left(n^2\right)$ elementary operations, and so we have fast algorithms, 
whereas the total computational cost of computing the eigenvalues or the singular values is $\mathcal{O}\left(n^3\right)$ elementary operations.
\end{rem}

\subsection{HRA computations with lattice path matrices}

If $\alpha,\beta>0$ and $\gamma\ge 0$, by Corollary \ref{corSTP}, the matrices $Lp_n(\alpha,\beta,\gamma)$ are STP and their bidiagonal decompositions can be computed to HRA,
and so the algebraic computations mentioned before can also be performed to HRA.

Let us consider the matrices 
$Lp_n(\sqrt{2},\sqrt{3},\sqrt{5})$ for $n=5,10,\ldots,50$.
First, we have computed the eigenvalues and the singular values of these matrices with Mathematica
using a precision of $100$ digits. We have also computed approximations to the eigenvalues of those matrices
in Matlab with \verb"eig" and also with \verb"TNEigenValues" using the
bidiagonal decomposition provided by \eqref{BDgeneral}. Then
we have computed the relative errors of the approximations obtained considering the eigenvalues
obtained with Mathematica as exact computations.

\begin{figure}
	\centering
    \includegraphics[width=1\textwidth]{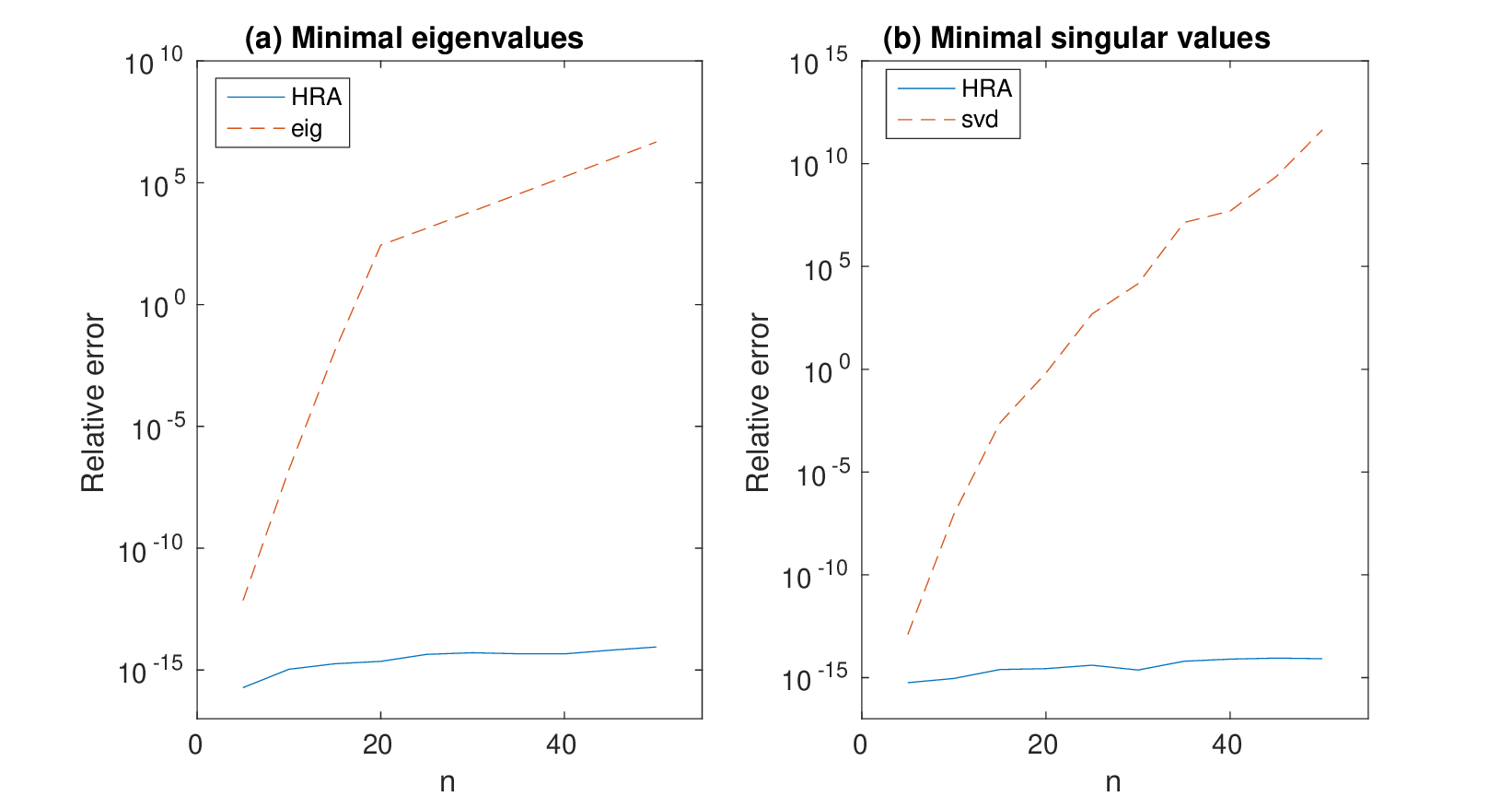}
  \caption{Relative error for the minimal eigenvalues and singular values of $Lp_n(\sqrt{2},\sqrt{3},\sqrt{5})$}
  \label{fig:ejemplo1}
\end{figure}

In Figure \ref{fig:ejemplo1} (a) we can see the relative error for the minimal eigenvalue of each matrix 
$Lp_n(\sqrt{2},\sqrt{3},\sqrt{5})$, $n=5,10,\ldots,50$, 
for both \verb"eig" and \verb"TNEigenValues". We can observe that Matlab function \verb"eig" does not 
provide an acceptable approximation of the minimal eigenvalue of the matrices $Lp_n(\sqrt{2},\sqrt{3},\sqrt{5})$
for $n\ge 15$ in contrast to the accurate approximations provided by the HRA computations of \verb"TNEigenValues".

We have also computed approximations to the singular values of the matrices 
$Lp_n(\sqrt{2},\sqrt{3},\sqrt{5})$, $n=5,10,\ldots,50$,
in Matlab with \verb"svd" and also with \verb"TNSingularValues". Then
we have computed the relative errors of the approximations obtained considering the singular values
obtained with Mathematica as exact computations.
In Figure \ref{fig:ejemplo1} (b) we can see the relative error for the minimal singular value of each matrix 
$Lp_n(\sqrt{2},\sqrt{3},\sqrt{5})$ 
for both \verb"svd" and \verb"TNSingularValues". As in the case of the eigenvalues, \verb"TNSingularValues"
provide very accurate approximations to the minimal singular values in contrast to the poor results
provided by \verb"svd".

\begin{center}
 \begin{figure*}[!htb]
 \includegraphics[keepaspectratio=true,width=1\textwidth]{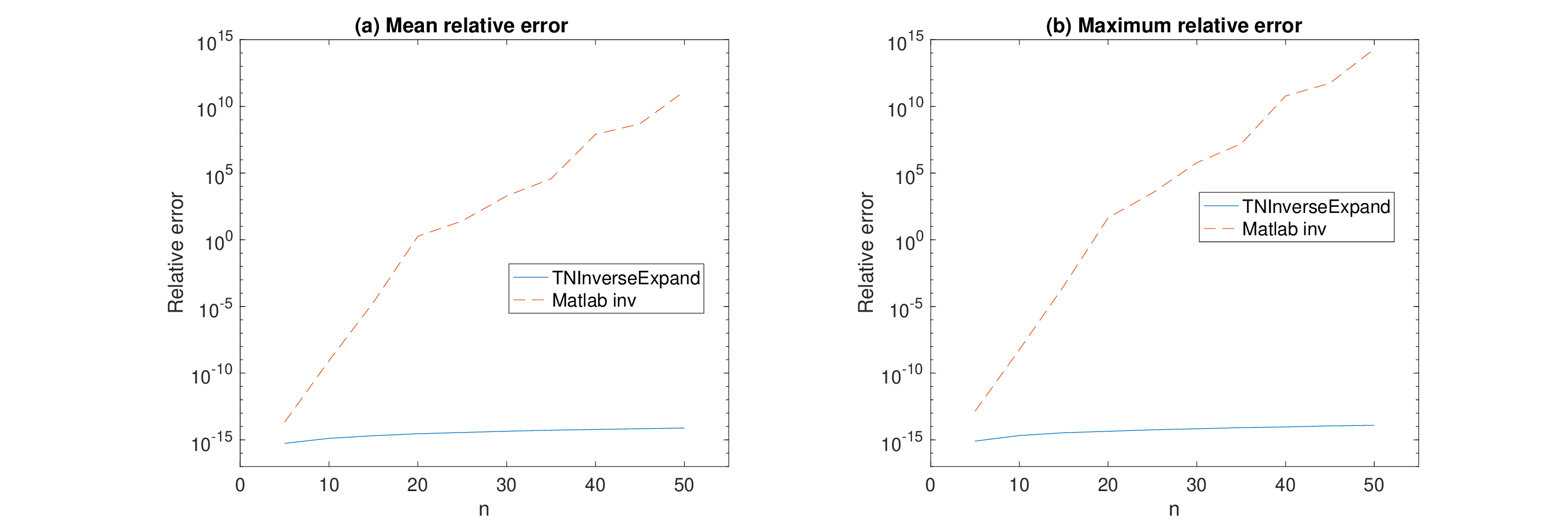}
  \caption{Relative errors for $Lp_n(\sqrt{2},\sqrt{3},\sqrt{5})^{-1}$, $n=5,10,\ldots,50$}\label{fig:inv1}
 \end{figure*}
\end{center}

\begin{center}
 \begin{figure*}[!h]
 \includegraphics[keepaspectratio=true,width=1\textwidth]{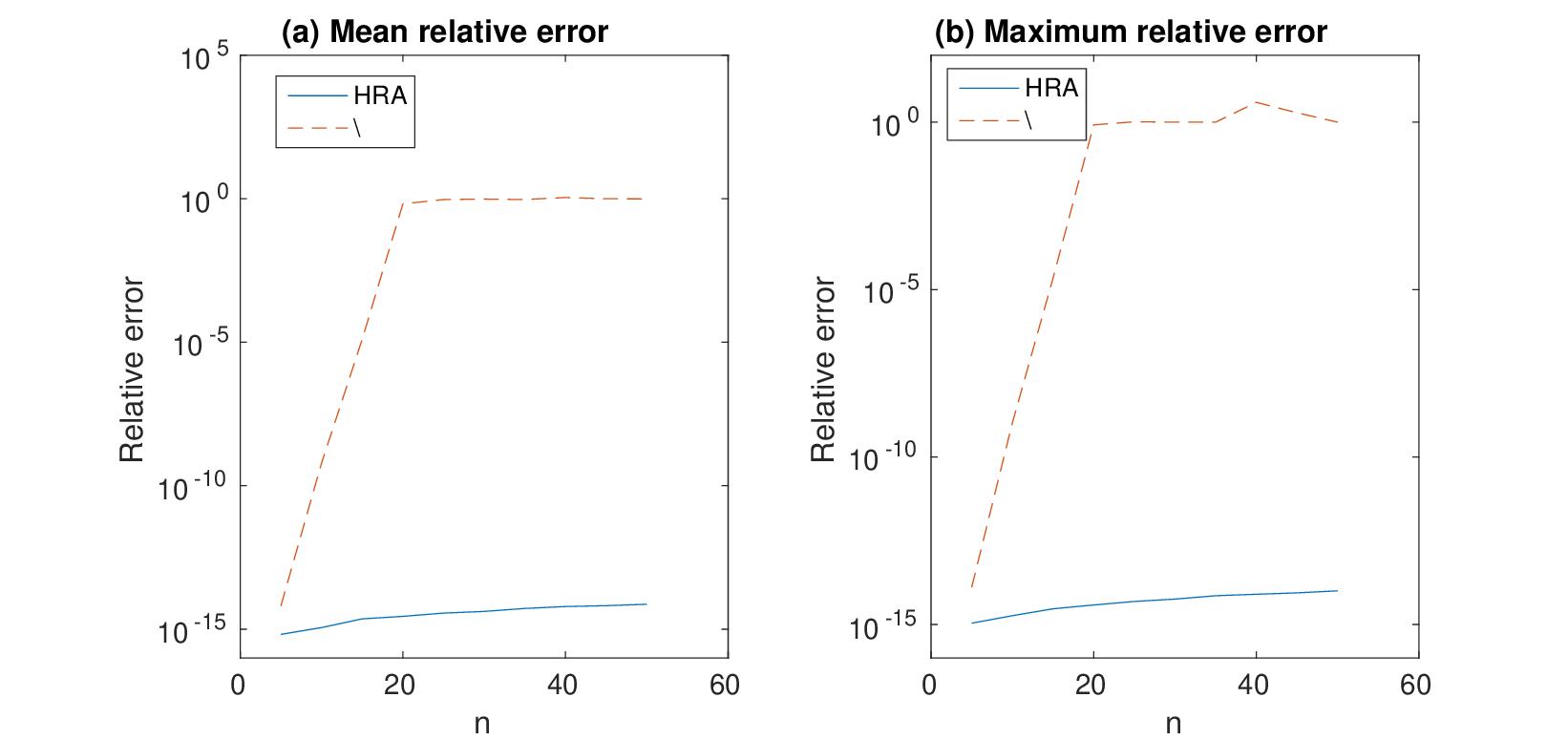}
  \caption{Relative errors for the systems $Lp_n(\sqrt{2},\sqrt{3},\sqrt{5})x=b_n$, $n=5,10,\ldots,50$}\label{fig:sist1}
 \end{figure*}
\end{center}

We have also computed with Matlab approximations to 
the inverses of the matrices $Lp_n(\sqrt{2},\sqrt{3},\sqrt{5})$, $n=5,10,\ldots,50$, with \verb"inv" and 
with \verb"TNInverseExpand" using the bidiagonal decomposition given
by \eqref{BDgeneral}. The 
inverses of these matrices have been computed with Mathematica using a precision of $100$ digits. 
Then we have computed the corresponding componentwise
relative errors. Finally we have obtained the mean and 
maximum componentwise relative 
errors. Figure \ref{fig:inv1} (a) shows the mean relative error and (b) 
shows the maximum relative error. We can also observe in this case 
that the results obtained with \verb"TNInverseExpand" are much more
accurate than those obtained with \verb"inv". In fact, the approximations
obtained with \verb"inv" are not acceptable for $n>15$.
Now we consider the linear systems $Lp_n(\sqrt{2},\sqrt{3},\sqrt{5})x=b_n$ for $n=5,10,\ldots,50$, where 
$b_n\in\mathbb{R}^n$ has the absolute value of its entries randomly
generated as integers in the interval 
$[1,1000]$, but with alternating signs. 
We have computed approximations to the solution $x$ of the linear systems with Matlab, the first one using
\verb"TNSolve" and the bidiagonal decomposition given
by \eqref{BDgeneral}, and the second one using the
Matlab command \verb"A\b".
By using Mathematica with a precision of $100$ digits we have computed the solution of the systems and then
we have computed the componentwise relative errors for the two approximations obtained with Matlab.
Then we have obtained the mean and maximum componentwise relative error. Figure \ref{fig:sist1} (a)
shows the mean relative error and (b) shows the maximum relative error.
Again, the results obtained with HRA algorithms are very accurate in contrast
to the results obtained with the usual Matlab command.

\begin{center}
 \begin{figure*}[!h]
 \includegraphics[keepaspectratio=true,width=1\textwidth]{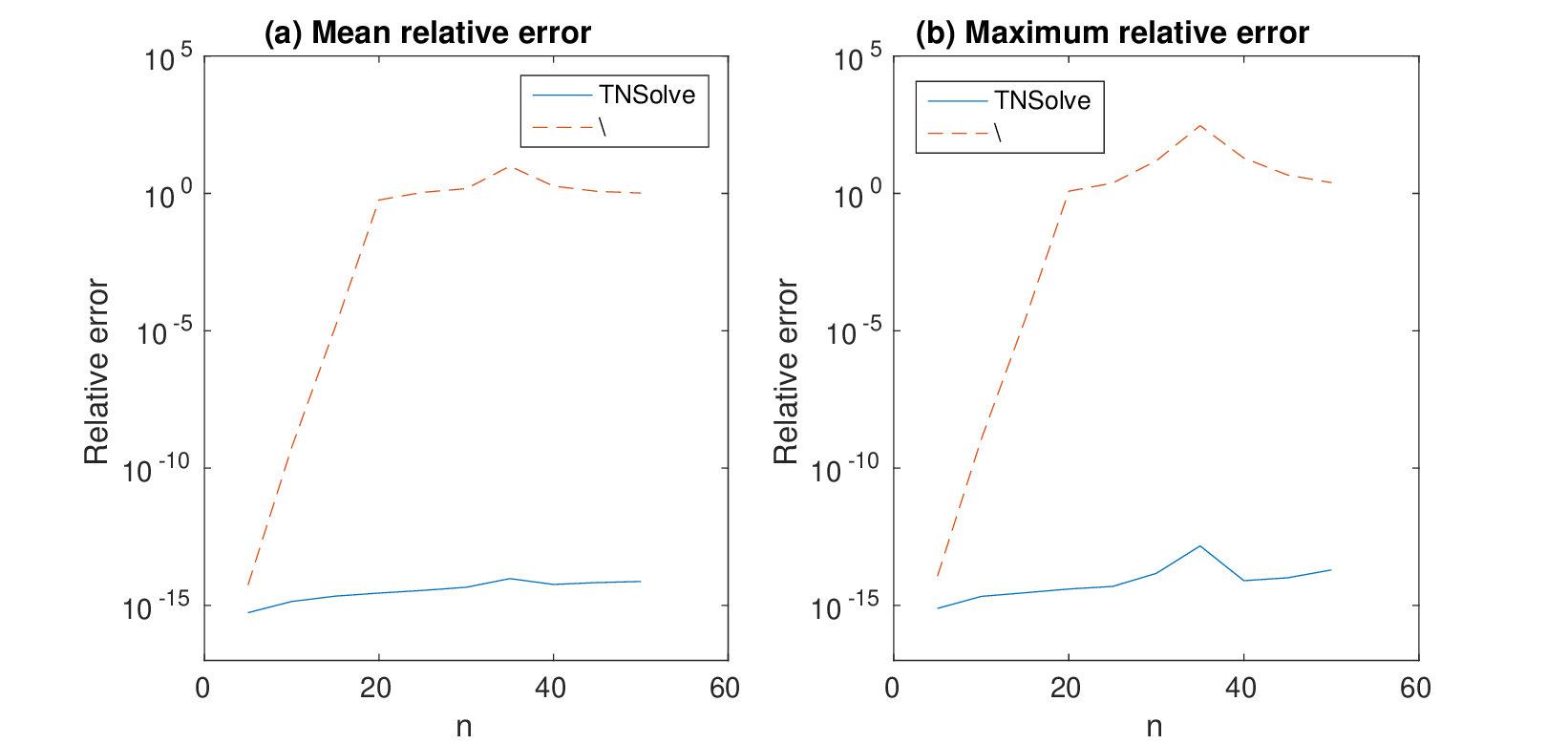}
  \caption{Relative errors for the systems $Lp_n(\sqrt{2},\sqrt{3},\sqrt{5})x=\tilde{b}_n$, $n=5,10,\ldots,50$}\label{fig:sist1.2}
 \end{figure*}
\end{center}

We also consider the linear systems $Lp_n(\sqrt{2},\sqrt{3},\sqrt{5})x=\tilde{b}_n$ for $n=5,10,\ldots,50$, 
where now $\tilde{b}_n\in\mathbb{R}^n$ has its entries randomly
generated as integers in the interval $[-1000,1000]$. 
 Figure \ref{fig:sist1.2} (a)
shows the mean relative error and (b) shows the maximum relative error.
In this case, we cannot  guarantee that 
the solution of the linear systems provided by \verb"TNSolve" can be computed to HRA.
However, the results obtained with \verb"TNSolve" are very accurate in contrast
to the results obtained with the usual Matlab command.

\subsection{Accurate computations for generalized triangular Pascal matrices}
Let us consider the lower triangular matrices $P_{n,1}[3/2]$ for $n=5,10,\ldots,50$, given by \eqref{gen3} with $x=3/2$ and $\lambda=1$. Unfortunately, by Corollary \ref{GPT.TP}, these matrices are not TP and we cannot assure that their bidiagonal decomposition can be computed to HRA. So we cannot guarantee that the algebraic computations mentioned above can be performed to HRA neither. Nevertheless,
let us also compare the numerical accuracy of \verb"TNSingularValues", \verb"TNInverseExpand"
and \verb"TNSolve" versus the usual Matlab commands \verb"svd", \verb"inv" and \verb"\", respectively.

\begin{figure}
	\centering
    \includegraphics[keepaspectratio=true,width=1\textwidth]{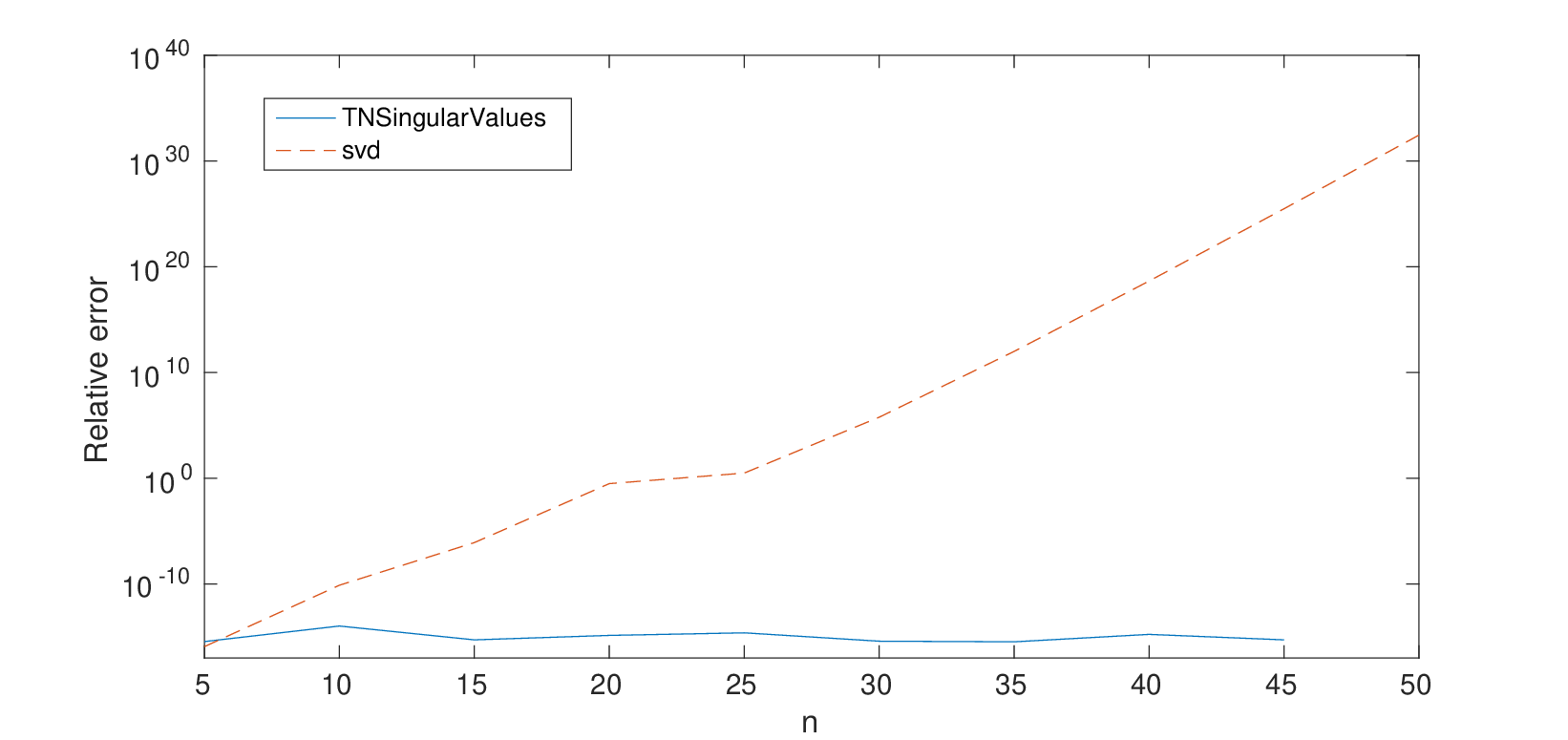}
  \caption{Relative error for the minimal singular values of $P_{n,1}[3/2]$}
  \label{fig:svd.ejemplo2}
\end{figure}

First, we have computed the singular values of these matrices with Mathematica
using a precision of $100$ digits. We have also computed approximations to the singular values of the
matrices $P_{n,1}[3/2]$ with Matlab function \verb"svd" and also with \verb"TNSingularValues" and
the corresponding $\mathcal{BD}(P_{n,1}[3/2])$ given in Theorem \ref{BD3}. Then
we have computed the relative errors of the approximations obtained considering the singular values
obtained with Mathematica as exact computations.
In Figure \ref{fig:svd.ejemplo2} we can see the relative error for the minimal singular value of each matrix 
for both \verb"svd" and \verb"TNSingularValues". 

\begin{center}
 \begin{figure*}[!htb]
 \includegraphics[keepaspectratio=true,width=1\textwidth]{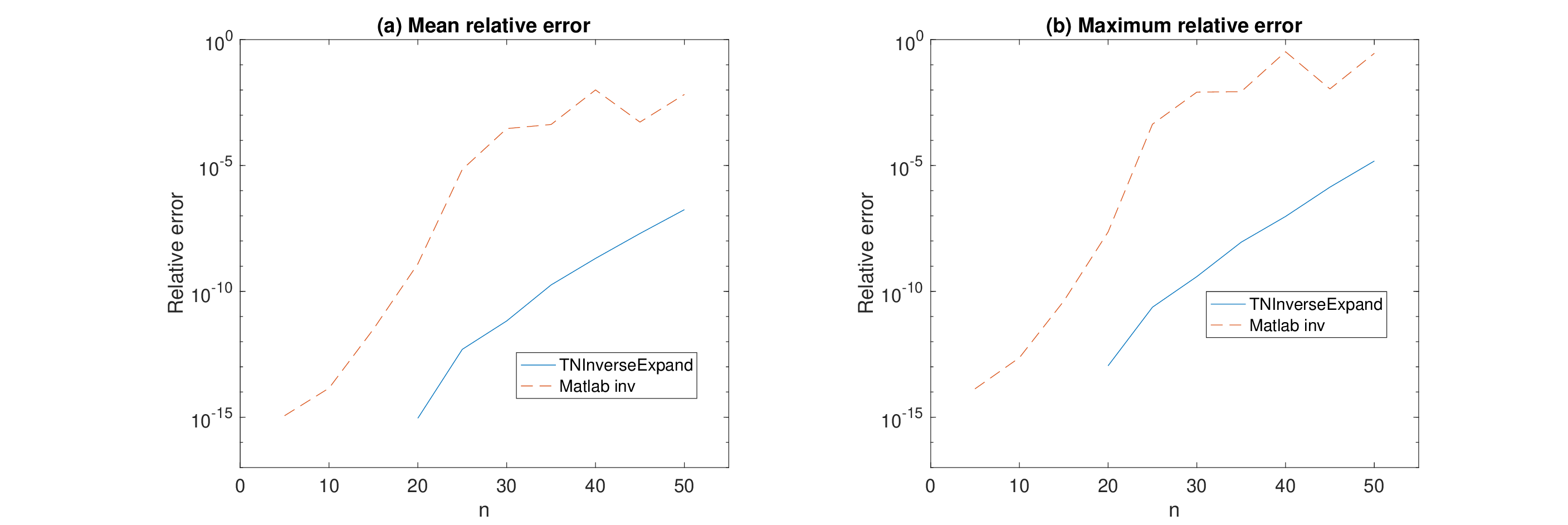}
  \caption{Relative errors for $P_{n,1}[3/2]^{-1}$, $n=5,10,\ldots,50$}\label{fig:inv1.ejemplo2}
 \end{figure*}
\end{center}

We have also computed with Matlab approximations to 
$P_{n,1}[3/2]^{-1}$, $n=5,10,\ldots,50$, with \verb"inv" and 
with \verb"TNInverseExpand" using $\mathcal{BD}(P_{n,1}[3/2])$. 
With Mathematica we have computed the 
inverse of these matrices with exact arithmetic. 
Then we have computed the corresponding componentwise
relative errors. Finally we have obtained the mean and 
maximum componentwise relative 
error. Figure \ref{fig:inv1.ejemplo2} (a) shows the mean relative error and (b) 
shows the maximum relative error. We can also observe in this case 
that the results obtained with \verb"TNInverseExpand" are much more
accurate than those obtained with \verb"inv".

\begin{center}
 \begin{figure*}[!h]
 \includegraphics[keepaspectratio=true,width=1\textwidth]{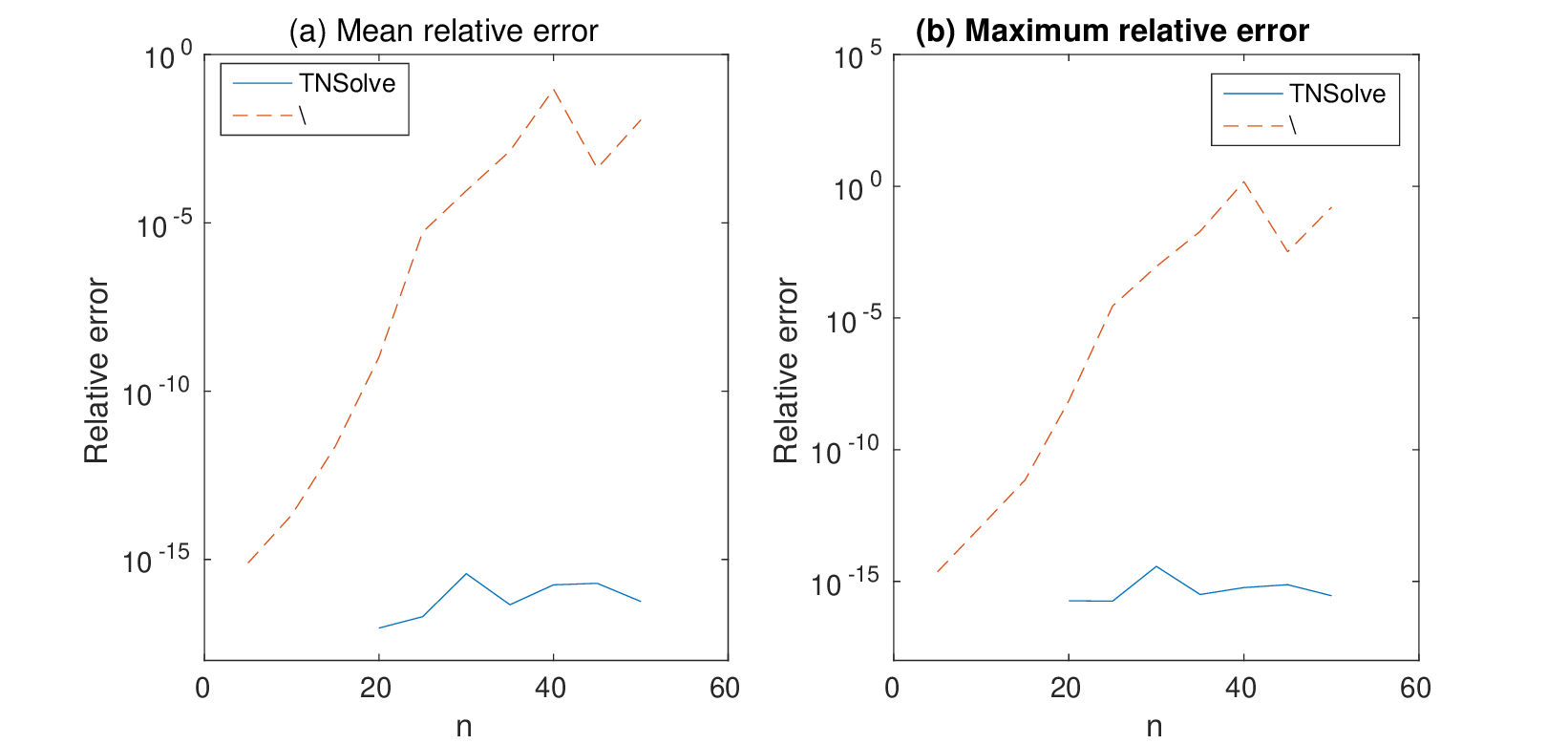}
  \caption{Relative errors for the systems $P_{n,1}[3/2]x=b_n$, $n=5,10,\ldots,50$}\label{fig:sist2}
 \end{figure*}
\end{center}

Finally we consider the linear systems $P_{n,1}[3/2]x=b_n$, $n=5,10,\ldots,50$, where 
$b_n\in\mathbb{R}^n$ has its entries randomly
generated as integers in the interval 
$[-1000,1000]$. We have computed approximations to the solution $x$ of the linear system with Matlab, the first one using
\verb"TNSolve" and $\mathcal{BD}(P_{n,1}[3/2])$, and the second one using the
Matlab command \verb"A\b".
By using Mathematica with exact arithmetic we have computed the exact solution of the systems and then
we have computed the componentwise relative errors for the two approximations obtained with Matlab.
Figure \ref{fig:sist2} (a)
shows the mean relative error and (b) shows the maximum relative error.
Again, the results obtained with \verb"TNSolve" are very accurate in contrast
to the results obtained with the usual Matlab command.

\section{Conclusions}\label{sec.num}
Pascal matrices and some generalizations considered in this paper arise in many applications, as commented in the introduction. It is well known that Pascal matrices and their generalizations are ill-conditioned (see \cite{a}). However, we show in this paper that we can compute with HRA
all their eigenvalues and all their singular values, and also the inverses of these matrices as well as the solutions of some linear systems. In fact, our 
 numerical experiments show that we can considerably improve the accuracy obtained with the usual Matlab commands. The crucial tool has been to obtain the bidiagonal decomposition of the generalized Pascal matrices with HRA and then apply the corresponding algorithms given in \cite{K1, K2, KoevSoft}. Let us also remark that, in spite of its much greater accuracy, the procedure presented in this paper has a computational cost similar to the usual algorithms used to solve these problems.




\begin{thebibliography}{}
 \bibitem{a} P. Alonso, J. Delgado, R. Gallego,  J.~M Pe\~na: Conditioning and accurate computations with Pascal matrices. J. Comput. Appl. Math. 252,  21--26 (2013).
\bibitem {A}  T. Ando: Totally positive matrices. Linear Algebra Appl. 90, 165--219 (1987)
\bibitem{Ref3} M. Bayat, H. Teimoori:
The linear algebra of the generalized Pascal functional matrix.
Linear Algebra Appl. 295, 81--89 (1999).
\bibitem{B} F. Brenti: Combinatorics and total positivity.
J. Combin. Theory Ser. A  71,  175--218 (1995). 
\bibitem{dop1} J. Delgado, H. Orera, J.~M. Pe\~na: Accurate computations with Laguerre matrices. Numer. Linear Algebra  Appl. 26: e2217 (10 pp.) (2019).
\bibitem{dop2} J. Delgado, H. Orera, J.~M. Pe\~na: Accurate algorithms for Bessel matrices. Journal of Scientific Computing 80, 1264--1278 (2019).
\bibitem{JS} J. Delgado, J.~M. Pe\~na: Fast and accurate algorithms for Jacobi-Stirling matrices. Appl. Math. Comput. 236, 253--259 (2014).
\bibitem{qBer} J. Delgado, J.~M. Pe\~na: Accurate computations with collocation matrices of q-Bernstein polynomials. SIAM J. Matrix Anal. Appl. 36, 880--893 (2015).
\bibitem{AN} J. Demmel,  I. Dumitriu,  O. Holtz,  P. Koev:
	Accurate and efficient expression evaluation and linear algebra. Acta Numer.  17,
	87--145 (2008).
\bibitem{DK} J. Demmel,  P. Koev: The accurate and efficient solution of a totally positive generalized Vandermonde linear system. SIAM J. Matrix Anal. Appl. 27 
142--152 (2005).
\bibitem{GM} M. Gasca, C.A. Micchelli (eds.): Total positivity and its applications. Mathematics and its Applications. Kluwer Acad. Publ. Dordrecht (1996).
\bibitem{TPandNE} M. Gasca,  J.~M. Pe\~na: Total positivity and Neville Elimination.  
Linear Algebra Appl. 165, 25--44 (1992)
\bibitem {fact} M. Gasca,  J.~M. Pe\~na: On factorizations of totally positive
matrices. In: {\it Total Positivity and Its Applications} (M.\,Gasca and C.A.\,Micchelli, Ed.),
Kluwer Academic Publishers, Dordrecht, The Netherlands (1996), pp. 109--130
\bibitem{Ref5} I.-P. Kim: LDU decomposition of an extension matrix of the Pascal matrix.
Linear Algebra Appl. 434, 2187–2196 (2011)
\bibitem{K1}  P. Koev: Accurate eigenvalues and SVDs of totally nonnegative matrices. SIAM J Matrix Anal. Appl. 27, 1--23 (2005)
\bibitem{K2}  P. Koev: Accurate computations with totally nonnegative matrices. SIAM J Matrix Anal. Appl. 29, 731--751 (2007)
\bibitem{KoevSoft}  P. Koev: \url{http://www.math.sjsu.edu/~koev/software/TNTool.html}. [Accessed 11 February 2020]
\bibitem{X} Lv, X.-G., Huang, T.-Z., Ren, Z.-G.: A new algorithm for linear systems of the Pascal type. J. Comput. Appl. Math. 225, 309--315 (2009).
\bibitem{SBV}  A. Marco, J.-J. Mart\'inez: Accurate computations with Said-Ball-Vandermonde matrices. Linear Algebra Appl. 432, 2894--2908 (2010).
\bibitem{BV2}  A. Marco, J.-J. Mart\'inez: Accurate computations with totally positive Bernstein-Vandermonde matrices. Electron. J. Linear Algebra 26, 357--380 (2013).
\bibitem{MM.Inv} A. Marco, J.-J. Mart\'inez J. J.:
Accurate computation of the Moore-Penrose inverse of strcitly
totally positive matrices. J. Comput. Appl. Math. 350, 299--308 (2019).
\bibitem{p} A. Pinkus: Totally positive matrices. Tracts in Mathematics; 181.  Cambridge University Press, Cambridge, UK (2010).
\bibitem{Ref1} Z. Zhang: The linear algebra of the generalized Pascal matrix. Linear Algebra Appl. 250, 51--60 (1997)
\bibitem{Ref2} Z. Zhang,  M. Liu: An extension of the generalized Pascal matrix and its algebraic properties. 
Linear Algebra Appl. 271, 169--177 (1998)

\end{thebibliography}
\end{document}